\documentclass[10pt]{amsart}

\usepackage{amsmath,amsthm,amsfonts,amscd,amssymb}
\usepackage{hyperref,wasysym}
\usepackage{verbatim}
\usepackage{amssymb}
\usepackage{graphicx}

\newtheorem{thm}{Theorem}[section]
\newtheorem{cor}[thm]{Corollary}
\newtheorem{lem}[thm]{Lemma}
\newtheorem{prop}[thm]{Proposition}

\theoremstyle{definition}

\theoremstyle{remark}

\begin{document}

\title{On lattice illumination of smooth convex bodies}

\author{Lenny Fukshansky}
\address{Department of Mathematics, 850 Columbia Avenue, Claremont McKenna College, Claremont, CA 91711}
\email{lenny@cmc.edu}

\subjclass[2020]{Primary: 52C07, 11H06, 52A05}
\keywords{illumination problem, lattice, orthogonality defect}

\begin{abstract} 
The illumination conjecture is a classical open problem in convex and discrete geometry, asserting that every compact convex body~$K$ in $\mathbb R^n$ can be illuminated by a set of no more than $2^n$ points. If $K$ has smooth boundary, it is known that $n+1$ points are necessary and sufficient. We consider an effective variant of the illumination problem for bodies with smooth boundary, where the illuminating set is restricted to points of a lattice and prove the existence of such a set close to $K$ with an explicit bound on the maximal distance. We produce improved bounds on this distance for certain classes of lattices, exhibiting additional symmetry or near-orthogonality properties. Our approach is based on the geometry of numbers.
\end{abstract}

\maketitle

\def\A{{\mathcal A}}
\def\B{{\mathcal B}}
\def\C{{\mathcal C}}
\def\D{{\mathcal D}}
\def\F{{\mathcal F}}
\def\x{{\mathcal H}}
\def\I{{\mathcal I}}
\def\J{{\mathcal J}}
\def\K{{\mathcal K}}
\def\L{{\mathcal L}}
\def\M{{\mathcal M}}
\def\N{{\mathcal N}}
\def\O{{\mathcal O}}
\def\R{{\mathcal R}}
\def\s{{\mathcal S}}
\def\V{{\mathcal V}}
\def\W{{\mathcal W}}
\def\X{{\mathcal X}}
\def\Y{{\mathcal Y}}
\def\H{{\mathcal H}}
\def\Z{{\mathcal Z}}
\def\OO{{\mathcal O}}
\def\BB{{\mathbb B}}
\def\cee{{\mathbb C}}
\def\EE{{\mathbb E}}
\def\Nn{{\mathbb N}}
\def\pee{{\mathbb P}}
\def\que{{\mathbb Q}}
\def\real{{\mathbb R}}
\def\zed{{\mathbb Z}}
\def\hyp{{\mathbb H}}
\def\aa{{\mathfrak a}}
\def\HH{{\mathfrak H}}
\def\qbar{{\overline{\mathbb Q}}}
\def\eps{{\varepsilon}}
\def\ahat{{\hat \alpha}}
\def\bhat{{\hat \beta}}
\def\gt{{\tilde \gamma}}
\def\h{{\tfrac12}}
\def\be{{\boldsymbol e}}
\def\bei{{\boldsymbol e_i}}
\def\bff{{\boldsymbol f}}
\def\ba{{\boldsymbol a}}
\def\bb{{\boldsymbol b}}
\def\bc{{\boldsymbol c}}
\def\bm{{\boldsymbol m}}
\def\bn{{\boldsymbol n}}
\def\bk{{\boldsymbol k}}
\def\bi{{\boldsymbol i}}
\def\bl{{\boldsymbol l}}
\def\bq{{\boldsymbol q}}
\def\bu{{\boldsymbol u}}
\def\bt{{\boldsymbol t}}
\def\bs{{\boldsymbol s}}
\def\bv{{\boldsymbol v}}
\def\bw{{\boldsymbol w}}
\def\bx{{\boldsymbol x}}
\def\bX{{\boldsymbol X}}
\def\bz{{\boldsymbol z}}
\def\bwy{{\boldsymbol y}}
\def\bY{{\boldsymbol Y}}
\def\bL{{\boldsymbol L}}
\def\baa{{\boldsymbol\alpha}}
\def\bbb{{\boldsymbol\beta}}
\def\bet{{\boldsymbol\eta}}
\def\bxi{{\boldsymbol\xi}}
\def\bo{{\boldsymbol 0}}
\def\bol{{\boldkey 1}_L}
\def\ep{\varepsilon}
\def\p{\boldsymbol\varphi}
\def\q{\boldsymbol\psi}
\def\rank{\operatorname{rank}}
\def\aut{\operatorname{Aut}}
\def\lcm{\operatorname{lcm}}
\def\sgn{\operatorname{sgn}}
\def\spn{\operatorname{span}}
\def\md{\operatorname{mod}}
\def\Norm{\operatorname{Norm}}
\def\dim{\operatorname{dim}}
\def\det{\operatorname{det}}
\def\Vol{\operatorname{Vol}}
\def\rk{\operatorname{rk}}
\def\Gal{\operatorname{Gal}}
\def\WR{\operatorname{WR}}
\def\WO{\operatorname{WO}}
\def\GL{\operatorname{GL}}
\def\pr{\operatorname{pr}}
\def\Tr{\operatorname{Tr}}
\def\dd{\partial}
\def\itt{\operatorname{int}}
\def\Ar{\operatorname{Area}}

\section{Introduction and statement of results}
\label{intro}

Let $K$ be a compact convex body in $\real^n$, write $\dd K$ for its boundary and $\itt(K)$ for its interior. Let $\bx \in \real^n$ and $\bwy \in \dd K$. We say that $\bx$ {\it illuminates} $\bwy$ if the line segment connecting $\bx$ and $\bwy$ does not intersect $\itt(K)$ but the line containing this line segment intersects $\itt(K)$. A collection of points $S \subset \real^n$ is said to illuminate $K$ if every point $\bwy \in \dd K$ is illuminated by some point $\bx \in S$. The illumination number of $K$ is defined as
$$I(K) := \min \{ |S| : S \text{ illuminates } K \},$$
where $|S|$ stands for the cardinality of the set $S$. The Illumination Conjecture then asserts that for any $n$-dimensional convex body $K$, $I(K) \leq 2^n$ with $I(K) = 2^n$ if and only if $K$ is an affine image of an $n$-cube (see~\cite{bikeev}, \cite{sriamorn} for more details). On the other hand, if $K$ has smooth boundary (i.e. there is a unique support hyperplane at each point $\bwy \in \dd K$), then it is well known that $I(K) = n+1$ (see \cite{hadwiger}). In fact, it is not difficult to see that in this case $K$ can be illuminated by vertices of a simplex containing $K$ in its interior.

In this note, we are interested in an effective version of this illumination problem. Let us write $\K_n$ for the set of all convex compact bodies with smooth boundary in $\real^n$ and assume $K \in \K_n$. For a finite set $S$ that illuminates $K$, define the {\it illumination distance}
$$d(S,K) := \max \{ \|\bx - \bwy\| : \bx \in S, \bwy \in K \},$$
where $\|\ \|$ stands for the Euclidean norm on $\real^n$. Define the diameter of $K$ as
$$D(K) := \max \{ \|\bx - \bwy\| : \bx,\bwy \in K \}.$$
Our first observation is that $K$ can be illuminated by a set of cardinality $n+1$ with bounded illumination distance. 

\begin{prop} \label{simplex} Let $K \in \K_n$, then for any $\eps > 0$ there exists $S \subset \real^n$ with $|S| = n+1$ that illuminates $K$ so that
$$d(S,K) \leq  \sqrt{\frac{n (n+1)}{2}} D(K) + \eps.$$
\end{prop} 

\noindent
We prove this proposition in Section~\ref{simplex_proof}. We inscribe $K$ into a ball and use Jung's inequality~\cite{jung} to bound the radius of this ball in terms of the diameter of $K$. We then illuminate the ball by the set of vertices of a regular simplex and show that the same set illuminates $K$.  
\smallskip

Now, let us consider a more delicate problem. Suppose $L \subset \real^n$ is a lattice of full rank. It is again not difficult to see that there exists a simplex with vertices on points of $L$ that illuminate $K$. The main goal of this note is to prove the existence of a set $S \subset L$ illuminating $K$ with bounded illumination distance. In what follows, we always identify any translated copy of $K$ with $K$. Hence, when we say that some set $S$ illuminates $K$, we mean that it illuminates some translated copy of $K$. We write $\det(L)$ for the determinant of the lattice $L$ and
$$\|L\| = \min \left\{ \|\bx\| : \bx \in L \setminus \{\bo\} \right\}$$
for the minimal norm of $L$.

\begin{thm} \label{main} Let $K \in \K_n$ and let $L \subset \real^n$ be a lattice of full rank. There exists $S \subset L$ with $|S| = n+1$ that illuminates $K$ so that
$$d(S,K) \leq 2 \left( \frac{4}{3} \right)^{\frac{n(n-1)}{2}} \left[ \left( \frac{4}{3} \right)^{\frac{n(n-1)}{2}} \frac{C_n D(K)}{\|L\|} + 1 \right] \frac{\det(L)}{\|L\|^{n-1}},$$
where $[\ ]$ stands for the integer part and the dimensional constant $C_n$ is given by
\begin{equation}
\label{Cn}
C_n = \frac{ n \sqrt{n} \left( n+2^{n-1} \right)}{\sqrt{2(n+1)}}.
\end{equation}
\end{thm}

\noindent
The set $S$ in the statement of this theorem is the set of vertices of a lattice simplex containing $K$ in its interior. To prove Theorem~\ref{main} in Section~\ref{lattice}, we construct such a simplex using an isoperimetric inequality and techniques from the geometry of numbers. Most importantly, we express the upper bound on the illumination distance in terms of the orthogonality defect of the illuminating set, a notion we define in Section~\ref{lattice}. We also explore our bound in more details for specific classes of lattices; all the notation discussed below is carefully reviewed in Section~\ref{lattice}.
\smallskip

The bound of Theorem~\ref{main} involves a dimensional constant, the unavoidable dependence on the diameter of $K$ and the dependence on the lattice $L$. We can minimize the dependence on $L$ in the case of well-rounded lattices; recall that $L$ is called well-rounded if it contains $n$ linearly independent vectors of minimal norm~$\|L\|$.

\begin{cor} \label{WR} With notation of Theorem~\ref{main}, assume that the lattice $L$ is well-rounded. Then
$$d(S,K) \leq 2 \left[ \frac{2^n C_n D(K)}{\omega_n \|L\|} + 1 \right] \|L\|,$$
where $C_n$ is as in~\eqref{Cn} and $\omega_n$ is the volume of a unit ball in $\real^n$.
\end{cor}

\noindent
Corollary~\ref{WR} suggests that a symmetric property like well-roundedness can potentially reduce illumination distance. Another property that can also be beneficial in this context is near-orthogonality; we review it in details in Section~\ref{lattice}.

\begin{cor} \label{near_orth} With notation of Theorem~\ref{main}, assume that the lattice $L$ is nearly orthogonal. Then
$$d(S,K) \leq 2 \left( \frac{4}{3} \right)^{\frac{n-1}{2}} \left[ \left( \frac{4}{3} \right)^{\frac{n-1}{2}} \frac{ C_n D(K) }{\|L\|} + 1 \right] \frac{\det(L)}{\|L\|^{n-1}},$$
where $C_n$ is as in~\eqref{Cn}. If in addition $L$ is well-rounded, then
$$d(S,K) \leq 2 \left[ \left( \frac{4}{3} \right)^{\frac{n-1}{2}} \frac{ C_n D(K) }{\|L\|} + 1 \right] \|L\|.$$
\end{cor}

\noindent
Finally, if a lattice is virtually rectangular, i.e., contains a finite-index sublattice with an orthogonal basis, we can obtain a bound with a smaller dimensional constant at the expense of a higher power on the determinant of this lattice.

\begin{cor} \label{virt_rect} With notation of Theorem~\ref{main}, assume that the lattice $L$ is virtually rectangular. Then it is isometric to a lattice of the form $CB\zed^n$, where $C$ is a nonsingular diagonal matrix and $B$ is a nonsingular integer matrix with relatively prime entries in each row. In this case,
$$d(S,K) \leq 2 \left[ \frac{ C_n D(K) }{\|L\|} + 1 \right] \frac{\det(L)^n}{|c|^{n-1} \|L\|^{n-1}},$$
where $C_n$ is as in~\eqref{Cn}.
\end{cor}

\noindent
We are now ready to proceed.

\bigskip

\section{Illumination distance}
\label{simplex_proof}

In this section, we prove Proposition~\ref{simplex}. Let $K \in \K_n$. Our first observation asserts that, if a set $S$ illuminates a ball containing $K$, then it illuminates $K$. While this is well understood, we include it here for the purposes of self-containment.

\begin{lem} \label{ball1} Let $B$ be any ball in $\real^n$ so that $K \subseteq B$. Let $S \subset \real^n$ be a finite set illuminating~$B$. Then $S$ illuminates~$K$.
\end{lem}

\begin{figure}
\centering
\includegraphics[scale=0.4]{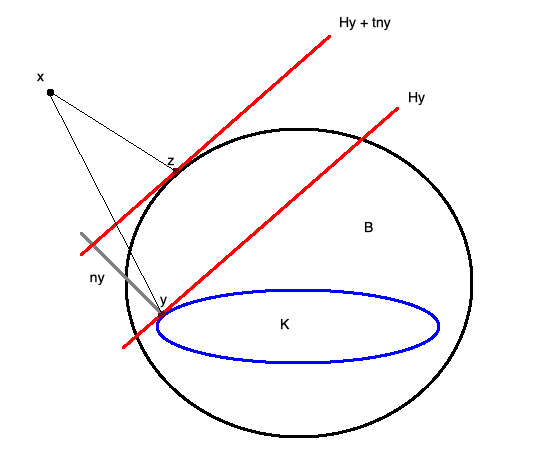}
\caption{Illustration of the construction in the proof of Lemma~\ref{ball1}.}\label{ball_fig}
\end{figure}

\proof
Let $\bwy \in \dd K$. We want to show that there exists some $\bx \in S$ such that $\bx$ illuminates $\bwy$. Let $\bn_{\bwy}$ be the outward-pointing unit normal vector to the support hyperplane $H_{\bwy}$ at $\bwy$. Let $\hyp^1_{\bwy}$ and $\hyp^2_{\bwy}$ be two open half-spaces so that $\real^n = \hyp^1_{\bwy} \sqcup H_{\bwy} \sqcup \hyp^2_{\bwy}$ and $K \subset \hyp^1_{\bwy}$. Then some point $\bx \in \real^n$ illuminates $\bwy$ if and only if $\bx \in \hyp^2_{\bwy}$. Let $t \in \real_{>0}$ be such that $H_{\bwy}+t\bn_{\bwy}$ is a support hyperplane to $B$ at some point $\bz \in \dd B$ and let $\hyp^1_{\bz}$, $\hyp^2_{\bz}$ be the corresponding open half-spaces, as above, so that $B \subset \hyp^1_{\bz}$. Then $\hyp^2_{\bz} \subset \hyp^2_{\bwy}$. Since the set $S$ illuminates $B$, there must be some point $\bx \in S \cap \hyp^2_{\bz}$ such that $\bx$ illuminates $\bz$. But $\bx \in \hyp^2_{\bwy}$ and so it illuminates $\bwy$. This construction is illustrated by Figure~\ref{ball_fig}.
\endproof

Next, we can describe a class of sets illuminating a ball $B$ in $\real^n$.

\begin{lem} \label{ball2} Let $P$ be a convex compact polytope in $\real^n$ containing $B$ in its interior. Then its set of vertices $S$ illuminates $B$.
\end{lem}

\proof
Let $\bwy \in \dd B$ and let $H_{\bwy}$ be the supporting hyperplane at $\bwy$. Let $\hyp^1_{\bwy}$ and $\hyp^2_{\bwy}$ be the two corresponding open half-spaces so that $B \subset \hyp^1_{\bwy}$. Since $B$ is contained in the interior of $P$, $H_{\bwy}$ intersects the interior of $P$. Hence, at least one vertex $\bv$ of $P$ is contained in $\hyp^2_{\bwy}$, and so it illuminates $\bwy$. Since this is true for any point in $\dd B$, the set of vertices of $P$ illuminates $B$.
\endproof

\begin{lem} \label{ball3} Let $B$ be a ball of radius $R$ in $\real^n$. Then for every $\eps > 0$ there exists a set $S$ of $n+1$ points in $\real^n$ illuminating $B$ so that
$$d(S,B) \leq (n+1) R + \eps.$$
\end{lem}

\proof
Let $P_t$ be a regular simplex with side length $t$. Then its inradius is given by the formula
$$r(P_t) = \frac{t}{\sqrt{2n(n+1)}}$$
and its height is
$$h(P_t) = t \sqrt{\frac{n+1}{2n}}.$$
Notice that if $r(P_t) > R$, then $P_t$ contains a copy of the ball $B$ in its interior. Hence, Lemma~\ref{ball2} implies that the set $S_t$ of $n+1$ vertices of $P_t$ illuminates $B$. For this to hold, we need to have
$$t > R \sqrt{2n(n+1)}.$$
Notice that $h(P_t)$ is the maximal distance from a vertex of $P_t$ to a point in $B$. Then let $\eps > 0$ and take $t \leq R \sqrt{2n(n+1)} + \eps \sqrt{\frac{2n}{n+1}}$, then for the corresponding set of vertices $S_t$ of $P_t$, we have
$$d(S_t,B) = h(P_t) \leq (n+1) R + \eps.$$
Take $S$ to be $S_t$ for any such choice of $t$, and this completes the proof.
\endproof

\proof[Proof of Proposition~\ref{simplex}]
Let $\BB(K) \subset \real^n$ be the smallest ball containing $K$. By Jung's inequality (see~\cite{jung}, also~\cite{circumradius} for a more contemporary account), the radius of $\BB(K)$ is
\begin{equation}
\label{jung_ineq}
R(K) \leq \sqrt{\frac{n}{2(n+1)}} D(K).
\end{equation}
By Lemma~\ref{ball1}, any set $S$ illuminating $\BB(K)$ illuminates $K$. By Lemma~\ref{ball3}, for any $\eps > 0$ there exists such a set $S$ with $|S| = n+1$ and 
$$d(S,\BB(K)) \leq (n+1) R(K) + \eps  \leq  \sqrt{\frac{n (n+1)}{2}} D(K) + \eps.$$
Since $K \subseteq B(K)$, $d(S,K) \leq d(S,\BB(K))$ and the result follows.
\endproof

\bigskip

\section{Lattice illumination}
\label{lattice}

In this section we prove Theorem~\ref{main} and its corollaries. Let $L \subset \real^n$ be a lattice of full rank. Let $\A = \{ \ba_1,\dots,\ba_n \}$ be a collection of linearly independent vectors in $L$, ordered so that
$$0 < \|\ba_1\| \leq \dots \leq \|\ba_n\|,$$
and let us write $\|\A\| = \|\ba_1\|$ for the minimal norm of these vectors. Let $A = (\ba_1\ \dots\ \ba_n)$ be the $n \times n$ nonsingular matrix with columns being the vectors of $\A$ and let $\Delta = |\det(A)|$. Define the orthogonality defect of $\A$ to be 
\begin{equation}
\label{orth_def}
\delta(\A) = \frac{\prod_{i=1}^n \|\ba_i\|}{\Delta}.
\end{equation}
By Hadamard's inequality, $\delta(\A) \geq 1$ with equality if and only if the collection $\A$ is orthogonal.
\smallskip

Let $t$ be a positive integer and let $S_{\A}(t)$ be the simplex whose vertices are $\bo, t \ba_1,\dots,t\ba_n \in L$. Write $\Vol_n(S_{\A}(t))$ for the volume of $S_{\A}(t)$ and $\Ar_{n-1}(S_{\A}(t))$ for its surface area. An isoperimetric inequality for $r(t)$ (see, e.g., (9) of \cite{hansen}), the inradius of $t$ guarantees that
\begin{equation}
\label{inradius}
r(t) \geq \frac{\Vol_n(S_{\A}(t))}{\Ar_{n-1}(S_{\A}(t))}.
\end{equation}
Notice that 
\begin{equation}
\label{s_vol}
\Vol_n(S_{\A}(t)) = \frac{|\det(tA)|}{n!} = \frac{t^n \Delta}{n!}.
\end{equation}
For each $1 \leq j \leq n$, let us write $F_j$ for the $(n-1)$-dimensional volume of the face of $S_{\A}(t)$ with vertices $\bo$ and $t \ba_i$ for all $i \neq j$. Since each such face is an $(n-1)$-dimensional simplex, the Hadamard inequality provides
\begin{equation}
\label{Fj}
F_j \leq \frac{1}{(n-1)!} \prod_{i=1, i \neq j}^n \|t \ba_i\| = \frac{t^{n-1} \prod_{i=1}^n \|\ba_i\|}{(n-1)! \|\ba_j\|} \leq \frac{t^{n-1} \prod_{i=1}^n \|\ba_i\|}{(n-1)! \|\A\|} = \frac{t^{n-1} \delta(\A) \Delta}{(n-1)! \|\A\|}.
\end{equation}
The only remaining face of $S_{\A}(t)$ is the $(n-1)$-dimensional simplex with vertices $t \ba_i$ for all $1 \leq i \leq n$; we write $F_{n+1}$ for its $(n-1)$-dimensional volume, then again by the Hadamard inequality
\begin{eqnarray}
\label{Flast}
F_{n+1} & \leq & \frac{1}{(n-1)!} \prod_{i=2}^n \|t \ba_i - t \ba_1\| \leq \frac{t^{n-1}}{(n-1)!} \prod_{i=2}^n \left( \|\ba_i\| + \|\ba_1\| \right) \nonumber \\
& \leq & \frac{(2t)^{n-1}}{(n-1)!} \prod_{i=2}^n \|\ba_i\| = \frac{(2t)^{n-1} \delta(\A) \Delta}{(n-1)! \|\A\|}.
\end{eqnarray}
Combining~\eqref{Fj} and~\eqref{Flast}, we obtain a bound
\begin{equation}
\label{area}
\Ar_{n-1}(S_{\A}(t)) = \sum_{i=1}^{n+1} F_i \leq  \left( \frac{n+ 2^{n-1}}{(n-1)!} \right) \frac{t^{n-1} \delta(\A) \Delta}{\|\A\|},
\end{equation}
and combining this inequality with~\eqref{s_vol} and~\eqref{inradius}, we get
\begin{equation}
\label{inradius1}
r(t) \geq \frac{t \|\A\|}{n \left( n+ 2^{n-1} \right) \delta(\A)}.
\end{equation}

\begin{figure}
\centering
\includegraphics[scale=0.45]{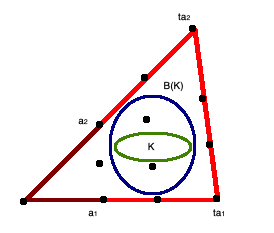}
\caption{Illustration of the construction in the proof of Theorem~\ref{main}.}\label{simplex_fig}
\end{figure}

Now, let $K \in \K_n$ and $\BB(K)$ be the smallest ball containing $K$. By Lemma~\ref{ball1}, any set illuminating $\BB(K)$ illuminates $K$. By Lemma~\ref{ball2}, if $\BB(K)$ is contained in the interior of $S_{\A}(t)$, then it is illuminated by the set of its vertices (see Figure~\ref{simplex_fig} for the illustration). In order for this to hold, we need to have the radius $R(K)$ for $\BB(K)$ to be smaller than $r(t)$, the inradius of $S_{\A}(t)$. To ensure this, we can take
$$\frac{t \|\A\|}{n \left( n+ 2^{n-1} \right) \delta(\A)} > R(K),$$
by~\eqref{inradius1}, i.e., $t > n \left( n+2^{n-1} \right) \frac{R(K)  \delta(\A) }{\|\A\|}$. For instance, we can take
\begin{equation}
\label{tstar}
t_* = \left[ n \left( n+2^{n-1} \right) \frac{R(K)  \delta(\A) }{\|\A\|} + 1 \right],
\end{equation}
and write $D(S_{\A}(t_*))$ for the diameter of $S_{\A}(t_*)$, which is given by
\begin{eqnarray*}
D(S_{\A}(t_*)) & = & \max \left\{ \|\bx - \bwy\| : \bx, \bwy \in S_{\A}(t_*) \right\} \\
& = & \max \left\{ \|t_* \ba_i - t_* \ba_j\| : \ba_i,\ba_j \in \A \right\} \leq 2t_* \|\ba_n\|,
\end{eqnarray*}
since $\| \ba_i - \ba_j \| \leq \|\ba_i\| + \|\ba_j\| \leq 2 \|\ba_n\|$. Further,
$$\|\ba_n\| = \frac{\delta(\A) \Delta}{\prod_{i=1}^{n-1} \|\ba_i\|} \leq  \frac{\delta(\A) \Delta}{\|\A\|^{n-1}}.$$
Now take $S$ to be the set of vertices of $S_{\A}(t_*)$ and observe that
$$d(S,K) \leq d(S,\BB(K)) \leq D(S_{\A}(t_*)) \leq 2t_* \frac{\delta(\A) \Delta}{\|\A\|^{n-1}}.$$
Combining this bound with~\eqref{tstar} and Jung's inequality~\eqref{jung_ineq}, we get
\begin{equation}
\label{ds_bnd}
d(S,K) \leq 2 \left[ n \left( n+2^{n-1} \right) \frac{ \sqrt{\frac{n}{2(n+1)}} D(K) \delta(\A) }{\|\A\|} + 1 \right] \frac{\delta(\A) \Delta}{\|\A\|^{n-1}}.
\end{equation}
Notice that in this bound $D(K)$ is certainly unavoidable. There is also the dimensional constant and the crucial dependence on $\delta(\A)$, $\Delta$ and $\|\A\|$ that we can attempt to minimize by an appropriate choice of the set $\A$.
\smallskip

\proof[Proof of Theorem~\ref{main}]
First notice that we can produce a bound dependent only on $L$, not on the choice of the set $\A$. Indeed, let $\A$ be an HKZ-reduced\footnote{Hermite-Korkin-Zolotarev reduction, see Section~2.9 of~\cite{martinet}} basis for $L$, then Hermite's inequality (see, e.g., \cite{martinet}, Section~2.2) guarantees that
\begin{equation}
\label{hermite}
\delta(\A) \leq \left( \frac{4}{3} \right)^{\frac{n(n-1)}{2}},
\end{equation}
although better, albeit more complicated, bounds are known (see~\cite{porter}). Since $\A$ is a basis for $L$,  $\Delta = \det(L)$. Further, $\|\A\| = \|L\|$, since an HKZ-reduced basis contains a shortest nonzero vector in the lattice. Then replacing $\delta(\A)$ with the bound~\eqref{hermite} in~\eqref{ds_bnd}, we obtain
\begin{equation}
\label{L_bnd}
d(S,K) \leq 2 \left( \frac{4}{3} \right)^{\frac{n(n-1)}{2}} \left[ \left( \frac{4}{3} \right)^{\frac{n(n-1)}{2}} \frac{C_n D(K)}{\|L\|} + 1 \right] \frac{\det(L)}{\|L\|^{n-1}},
\end{equation}
where $C_n$ is as in~\eqref{Cn}.

\endproof
\smallskip

We can do better for some special classes of lattices. We need some more notation. Given a full-rank lattice $L \subset \real^n$, let us write
$$\Sigma(L) = \left\{ \bx \in L : \|\bx\| = \|L\| \right\}$$
for its set of minimal vectors. $L$ is called {\it well-rounded} (WR) if $\Sigma(L)$ contains $n$ linearly independent vectors. 

\proof[Proof of Corollary~\ref{WR}]
If $L$ is WR, we can take $\A$ to be a set of linearly independent vectors contained in $\Sigma(L)$. Then
$$\frac{\Delta}{\|\A\|^n} = \frac{1}{\delta(\A)},$$
and Minkowski's Successive Minima Theorem (see, for instance, \cite[Section~9.1, Theorem 1]{gruber_lek}) gives a bound
$$\delta(\A) \leq \frac{2^n}{\omega_n}.$$
Combining these observations with~\eqref{ds_bnd}, we obtain
\begin{equation}
\label{L_bnd1}
d(S,K) \leq 2 \left[ \frac{2^n C_n D(K)}{\omega_n \|L\|} + 1 \right] \|L\|,
\end{equation}
where $C_n$ is as in~\eqref{Cn}.
\endproof
\smallskip

On the other hand, a lattice $L$ is called orthogonal if it possesses an orthogonal basis. More generally, let $\A = \{ \ba_1,\dots,\ba_n \}$ be an ordered basis for the lattice $L$, and define a sequence of angles $\theta_1,\dots,\theta_{n-1}$ as follows: each $\theta_i$ is the angle between $\ba_{i+1}$ and the subspace $\spn_{\real} \{ \ba_1,\dots,\ba_i \}$. It is then clear that each $\theta_i \in (0,\pi/2]$. If for each $1 \leq i \leq n-1$, $\theta_i \in [\theta,\pi/2]$ for some fixed $\theta \in (0, \pi/2]$ then $L$ is called (weakly) {\it $\theta$-orthogonal} ($L$ is $\theta$-orthogonal if the above condition holds for every reordering of the basis $\A$); we also refer to such lattices as {\it nearly orthogonal} for all $\theta \in [\pi/3,\pi/2]$. Nearly orthogonal lattices have been studied in~\cite{baraniuk} and~\cite{lf2}. 

\proof[Proof of Corollary~\ref{near_orth}]
Suppose $L$ is nearly orthogonal, then
$$\det(L) = \left( \prod_{i=1}^n \|\ba_i\| \right) \left( \prod_{i=1}^{n-1} \sin \theta_i \right) \geq \left( \prod_{i=1}^n \|\ba_i\| \right) (\sin \theta)^{n-1},$$
meaning that
$$\delta(\A) \leq \frac{1}{(\sin \theta)^{n-1}} \leq \left( \frac{2}{\sqrt{3}} \right)^{n-1} = \left( \frac{4}{3} \right)^{\frac{n-1}{2}}.$$
Hence, for nearly orthogonal lattice $L$ with this choice of $\A$,~\eqref{ds_bnd} gives
$$d(S,K) \leq 2 \left( \frac{4}{3} \right)^{\frac{n-1}{2}} \left[ \left( \frac{4}{3} \right)^{\frac{n-1}{2}} \frac{ C_n D(K) }{\|L\|} + 1 \right] \frac{\det(L)}{\|L\|^{n-1}}.$$
Suppose that $L$ is nearly orthogonal and WR. Then Theorem~1.1 of~\cite{lf2} guarantees that the nearly orthogonal basis consists of minimal vectors, and so
$$d(S,K) \leq 2 \left[ \left( \frac{4}{3} \right)^{\frac{n-1}{2}} \frac{ C_n D(K) }{\|L\|} + 1 \right] \|L\|.$$
\endproof
\smallskip

Two full-rank lattices $L, M \subset \real^n$ are called {\it isometric} if $M = U L$ for an $n \times n$ orthogonal matrix $U$. If a lattice $L$ contains a finite-index orthogonal sublattice, we call $L$ {\it virtually rectangular}. Virtually rectangular lattices were studied in~\cite{lf1}.

\proof[Proof of Corollary~\ref{virt_rect}]
Theorem~1.2 of~\cite{lf1} guarantees that $L$ is virtually rectangular if and only if it is isometric to a lattice of the form $\C B\zed^n$, where
$$\C = \begin{pmatrix} c_1 & \dots & 0 \\ \vdots & \ddots & \vdots \\ 0 & \dots & c_n \end{pmatrix}$$
with nonzero $c_1,\dots,c_n \in \real$ and $B$ a nonsingular integer matrix with relatively prime rows. Then
$$L = U \C L',$$
where $U$ is an $n \times n$ orthogonal matrix and $L' = B\zed^n \subseteq \zed^n$. Write $c = c_1 \cdots c_n$, then Theorem~1.3 of \cite{lf1} guarantees that there exists an orthogonal sublattice $\Lambda \subseteq L$ with
$$[L: \Lambda] = \left( \frac{\det(L)}{|c|} \right)^{n-1} = \left(  \det(L') \right)^{n-1}.$$
Then 
$$\det(\Lambda) = [L : \Lambda] \det(L) = \frac{\det(L)^n}{|c|^{n-1}}$$
and  $\|\Lambda\| \geq \|L\|$. Take $\A$ to be an orthogonal basis for $\Lambda$, then $\delta(\A)=1$ and~\eqref{ds_bnd} becomes
$$d(S,K) \leq 2 \left[ \frac{ C_n D(K) }{\|L\|} + 1 \right] \frac{\det(L)^n}{|c|^{n-1} \|L\|^{n-1}}.$$
\endproof

\bigskip

{\bf Acknowledgement:} I would like to thank Nitipon Moonwichit, from whom I learned about the illumination problem. I am also grateful to the anonymous referee for the many helpful remarks.

\bibliographystyle{plain}  

\end{document}